\documentclass[11pt,a4paper]{article}
\usepackage[T1]{fontenc}
\usepackage{lmodern}
\usepackage[margin=1in]{geometry}
\usepackage{amsmath,amssymb,amsthm}
\usepackage{microtype}
\usepackage{enumitem}
\usepackage{cite}
\usepackage[hidelinks]{hyperref}
\hypersetup{pdftitle={Understanding SVRG: Note IV on Variance Reduction}}
\newtheorem{assumption}{Assumption}
\newtheorem{lemma}{Lemma}
\newtheorem{theorem}{Theorem}
\newcommand{\E}{\mathbb{E}}
\newcommand{\R}{\mathbb{R}}
\newcommand{\norm}[1]{\left\|#1\right\|}
\newcommand{\ip}[2]{\left\langle #1,#2\right\rangle}
\newcommand{\eps}{\varepsilon}
\setlist[itemize]{leftmargin=1.4em,itemsep=2pt,topsep=4pt,parsep=0pt}
\title{Variance Reduction Note 4: The Convergence of SVRG}
\author{Alina Lodhi\\
xieyunqi27@gmail.com}
\date{}

\begin{document}
\maketitle
\vspace{-2em}
\begin{abstract}
Stochastic variance-reduced gradient (SVRG) uses a full gradient at a
snapshot point to correct inexpensive component gradients. This fourth
expository note on variance reduction explains that control variate and
gives a self-contained proof of linear convergence for the classical
strongly convex finite-sum setting. Each component is assumed convex
and smooth, while only their average must be strongly convex. A
second-moment bound links the estimator to the objective gaps at the
current iterate and the snapshot. Summing a distance recursion over
an epoch then yields a geometric decrease and an oracle complexity
of order $(n+L/\mu)\log(1/\tau)$ for objective accuracy $\tau$,
with the initial gap fixed. The result uses a random snapshot update
and counts all full-gradient computations. No new algorithm or rate
is claimed.
\end{abstract}

\section{A snapshot as a control variate}
Consider the finite-sum problem
\begin{equation}\label{eq:problem}
 \min_{x\in\R^d}F(x):=\frac1n\sum_{i=1}^n f_i(x).
\end{equation}
Sampling a single component gradient is inexpensive, but its noise can
persist at the minimizer: individual $\nabla f_i(x_\star)$ need not
vanish even when $\nabla F(x_\star)=0$. SVRG, introduced by Johnson
and Zhang \cite{svrg2013}, corrects this noise using a snapshot
$\widetilde x$ and its exact gradient. For a uniform random index $i$,
the estimator is
\begin{equation}\label{eq:estimator}
 v(x,\widetilde x;i)=\nabla f_i(x)-\nabla f_i(\widetilde x)
                          +\nabla F(\widetilde x).
\end{equation}
The subtracted and added terms have the same expectation. Consequently,
the estimator is unbiased for $\nabla F(x)$, while the shared index
cancels much of the component noise. At $x=\widetilde x$, it equals
the full gradient for every index.

This note develops the classical strongly convex analysis. Related
finite-sum methods include SAG, SAGA, and proximal SVRG
\cite{sag2012,saga2014,proxsvrg2014}. Recursive approaches such as
SARAH, SPIDER, and PAGE \cite{sarah2017,spider2018,page2021}
use a different form of correction. An overview of these method
families appears in \cite{gower2020}.

\section{Assumptions and the epoch algorithm}
\begin{assumption}\label{ass:convex}
Each $f_i:\R^d\to\R$ is convex, differentiable, and has an
$L$-Lipschitz gradient, with $L>0$. Their average $F$ is
$\mu$-strongly convex for some $0<\mu\leq L$: for all $x,y$,
\begin{equation}\label{eq:strong}
 F(y)\geq F(x)+\ip{\nabla F(x)}{y-x}
                    +\frac\mu2\norm{y-x}^2.
\end{equation}
Let $x_\star$ be its unique minimizer, $F_\star=F(x_\star)$,
and $h(x)=F(x)-F_\star$. Norms are Euclidean.
\end{assumption}
The components need not be strongly convex and need not share a
minimizer. We do not assume uniformly bounded stochastic-gradient
variance. The proof instead derives a bound that shrinks with
the objective gaps.

Fix an initial snapshot $\widetilde x_0$, a step size $\eta>0$, and
an integer epoch length $m\geq1$. At epoch $s=0,1,\ldots$:
\begin{enumerate}[leftmargin=1.6em,itemsep=2pt,topsep=4pt,parsep=0pt]
 \item Compute the full snapshot gradient
 $\widetilde g_s=n^{-1}\sum_{i=1}^n\nabla f_i(\widetilde x_s)$,
 and initialize $x_0=\widetilde x_s$.
 \item For $t=0,\ldots,m-1$, draw a fresh independent uniform
 index $i_t\in\{1,\ldots,n\}$ and update
 \begin{equation}\label{eq:update}
  v_t=\nabla f_{i_t}(x_t)-\nabla f_{i_t}(\widetilde x_s)
                 +\widetilde g_s,\qquad
  x_{t+1}=x_t-\eta v_t.
 \end{equation}
 \item Choose $J_s$ independently and uniformly from
 $\{0,\ldots,m-1\}$, and set $\widetilde x_{s+1}=x_{J_s}$.
\end{enumerate}
The inner index $t$ restarts at each epoch. This random snapshot
choice is the variant analyzed below; the proof averages the function
values of precisely these $m$ iterates. An implementation may draw
$J_s$ in advance and retain its selected iterate without storing
the whole trajectory. The cost per epoch is at most $n+2m$
component-gradient evaluations, using $O(d)$ auxiliary vector storage
when snapshot component gradients are recomputed as needed.

\section{A short linear-convergence proof}
The proof first controls component-gradient differences at the solution,
then bounds the SVRG estimator, and finally sums a distance recursion
over one epoch.

\begin{lemma}[Gradient differences and objective gaps]\label{lem:growth}
Under Assumption~\ref{ass:convex}, for every $x$,
\begin{equation}\label{eq:growth}
 \frac1n\sum_{i=1}^n
 \norm{\nabla f_i(x)-\nabla f_i(x_\star)}^2\leq2Lh(x).
\end{equation}
\end{lemma}
\begin{proof}
Define
$\psi_i(u)=f_i(u)-f_i(x_\star)
-\ip{\nabla f_i(x_\star)}{u-x_\star}$.
Convexity gives $\psi_i(u)\geq0$, and $\psi_i$ is $L$-smooth.
Apply its descent inequality at $u=x$ with
$y=x-\nabla\psi_i(x)/L$:
\[
 0\leq\psi_i(y)\leq\psi_i(x)
                    -\frac1{2L}\norm{\nabla\psi_i(x)}^2.
\]
Thus $\norm{\nabla f_i(x)-\nabla f_i(x_\star)}^2\leq2L\psi_i(x)$.
Averaging over $i$ gives \eqref{eq:growth}, because
$n^{-1}\sum_i\nabla f_i(x_\star)=\nabla F(x_\star)=0$.
\end{proof}

\begin{lemma}[Unbiasedness and second moment]\label{lem:moment}
For fixed $x,\widetilde x$ and a uniform index $i$, the estimator
in \eqref{eq:estimator} satisfies
\begin{equation}\label{eq:moment}
 \E_i[v]=\nabla F(x),\qquad
 \E_i\norm{v}^2\leq4L[h(x)+h(\widetilde x)].
\end{equation}
\end{lemma}
\begin{proof}
Unbiasedness follows by averaging \eqref{eq:estimator}. Set
$a_i=\nabla f_i(x)-\nabla f_i(x_\star)$ and
$c_i=\nabla f_i(\widetilde x)-\nabla f_i(x_\star)$.
Then $v=a_i-(c_i-\E_i c_i)$. Hence
\begin{align*}
 \E_i\norm{v}^2
 &\leq2\E_i\norm{a_i}^2+2\E_i\norm{c_i-\E_i c_i}^2\\
 &\leq2\E_i\norm{a_i}^2+2\E_i\norm{c_i}^2
 \leq4L[h(x)+h(\widetilde x)],
\end{align*}
using the variance identity and Lemma~\ref{lem:growth}.
\end{proof}

\begin{theorem}[Contraction per epoch]\label{thm:linear}
Suppose Assumption~\ref{ass:convex} holds. If $0<\eta<1/(2L)$
and $m$ are chosen so that
\begin{equation}\label{eq:rho}
 \rho:=\frac1{\mu\eta m(1-2L\eta)}
              +\frac{2L\eta}{1-2L\eta}<1,
\end{equation}
then after $S$ epochs,
\begin{equation}\label{eq:linear}
 \E h(\widetilde x_S)\leq\rho^S h(\widetilde x_0).
\end{equation}
In particular, $\eta=1/(10L)$ and
$m=\lceil50L/\mu\rceil$ ensure $\rho\leq1/2$.
\end{theorem}
\begin{proof}
Condition on the start of an epoch and denote its snapshot by
$\widetilde x$. All expectations in this paragraph are conditional
on that snapshot. Condition further on the history before drawing $i_t$.
Unbiasedness, convexity of $F$, and Lemma~\ref{lem:moment} imply
\begin{align*}
 \E_{i_t}\norm{x_{t+1}-x_\star}^2
 &=\norm{x_t-x_\star}^2
 -2\eta\ip{\nabla F(x_t)}{x_t-x_\star}
 +\eta^2\E_{i_t}\norm{v_t}^2\\
 &\leq\norm{x_t-x_\star}^2
 -2\eta(1-2L\eta)h(x_t)+4L\eta^2h(\widetilde x).
\end{align*}
Take expectations over the inner history and sum for $t=0,\ldots,m-1$.
The squared distances telescope, and the final one is nonnegative:
\[
 2\eta(1-2L\eta)\sum_{t=0}^{m-1}\E h(x_t)
 \leq\norm{\widetilde x-x_\star}^2
                       +4L\eta^2m h(\widetilde x).
\]
Strong convexity gives
$\norm{\widetilde x-x_\star}^2\leq2h(\widetilde x)/\mu$.
Divide by $2\eta(1-2L\eta)m$ and use the independent uniform
snapshot choice to obtain
$\E[h(\widetilde x_{s+1})\mid\widetilde x_s]
\leq\rho h(\widetilde x_s)$.
The tower property and induction prove \eqref{eq:linear}.
For the stated parameters, $1-2L\eta=4/5$, and the two
terms in \eqref{eq:rho} are at most $1/4$ and $1/4$.
\end{proof}

\paragraph{Oracle complexity and accuracy.}
Let $\kappa=L/\mu$ and $\Delta_0=h(\widetilde x_0)$.
With the parameters in Theorem~\ref{thm:linear}, an objective target
$0<\tau<\Delta_0$ is reached after
\begin{equation}\label{eq:epochs}
 S=\left\lceil\log_2\frac{\Delta_0}{\tau}\right\rceil
 \quad\Longrightarrow\quad
 \E h(\widetilde x_S)\leq\tau.
\end{equation}
The total component-gradient count is at most
\begin{equation}\label{eq:complexity}
 N\leq S(n+2\lceil50\kappa\rceil)
 =O\!\left((n+\kappa)
                \left[1+\log\frac{\Delta_0}{\tau}\right]\right).
\end{equation}
If $\Delta_0\leq\tau$, no epochs are needed for this theoretical
guarantee. Choosing $S$ in advance uses $\Delta_0$ or an upper bound.

The objective and stationarity criteria can also be related explicitly.
Applying smoothness at $x-\nabla F(x)/L$ and using the global lower
bound $F_\star$ gives
\begin{equation}\label{eq:gradient-target}
 \norm{\nabla F(x)}^2\leq2Lh(x).
\end{equation}
Thus taking $\tau=\eps^2/(2L)$ ensures
$\E\norm{\nabla F(\widetilde x_S)}^2\leq\eps^2$ and, by
Cauchy--Schwarz, $\E\norm{\nabla F(\widetilde x_S)}\leq\eps$.
Here $\tau$ measures objective error, while $\eps$ measures gradient
norm; the distinction matters when comparing rates across methods.

\paragraph{Interpretation and scope.}
The second moment in \eqref{eq:moment} decreases as both the inner
iterate and the snapshot approach the solution. This permits a fixed
step size and geometric improvement without assuming that every
component gradient vanishes at $x_\star$. The proof requires convex
components, a strongly convex average, and the specified random
snapshot rule. It does not establish the same rate for the last inner
iterate or a nonconvex objective. Nonconvex SVRG analyses
\cite{reddi2016,allenzhu2016} and loopless variants
\cite{loopless2019,looplesssampling2021} use different arguments.
The exact snapshot gradient is also essential to this proof;
approximating it introduces another error term.

\clearpage
\section{Recommended reading}

\paragraph{Variance-reduction methods.}
Basic guarantees \cite{ghadimi2013,gower2020,unifiedsgd2020}; gradient memories and control
variates \cite{sag2012,sdca2013,svrg2013,s2gd2013,saga2014,finito2014,proxsvrg2014}; sampling and acceleration \cite{catalyst2015,katyusha2016,svrgpp2016,katyushax2018,sagasampling2019,jacsketch2018,svrgsd2018,loopless2019,looplesssampling2021}; nonconvex and recursive
estimators \cite{reddi2016,allenzhu2016,natasha2017,sarah2017,isarah2018,scsg2017,spider2018,snvrg2018,spiderboost2019,proxsarah2020,zerosarah2021,page2021,nonconvexsampling2019}; momentum and adaptive schedules \cite{storm2019,hybrid2019,stormplus2021,superadam2021,metastorm2022,adaspider2022,jiang_adaptive2024,adasvrg2021,almosttune2020,adaptiveextra2022,alfpg2026};
second-order extensions \cite{ssrgd2019,natasha22018,svrc2019,svrcnewton2018}; lower bounds \cite{lower2023}.

\paragraph{Compositional optimization.}
Two-level  approximation \cite{scgd2014,asc2017,nasa2020,nearlyeasy2020}; finite-sum and recursive
estimators \cite{lian2017,compavg2018,civr2019,sarahcomp2019}; multilevel tracking \cite{nested2019,jiang_multilevel2022,levelindependent2020,jiang_revisiting2025};
coupled and partially observed blocks \cite{fccotheory2022,jiang_msvr2022,weakcomp2023,jiang_fcco2026,momentumcomp2025}; projection-free
methods \cite{vrfw2016,projectioncomp2022,jiang_projection2024,scgpp2019,jiang_projection2026}; statistical-model applications and federated
extensions \cite{jiang_models2023,fednest2022,jiang_models2026,localcomp2022}.

\paragraph{Sign-based optimization.}
Quantization, signSGD, majority voting, and error feedback \cite{qsgd2017,terngrad2017,signsgd2018,majority2019,errorfeedback2019};
variance-reduced sign updates and Lion \cite{jiang_sign2024,lion2023,jiang_lion2026}; compressed gradient
tracking \cite{diana2019,marina2021,jiang_distributed2026,ef212021,dasha2022}. Related distributed control-variate methods include
SCAFFOLD, FedPAGE, decentralized hybrid estimators, and ProxSkip \cite{scaffold2020,fedpage2021,hybriddecentral2021,proxskip2022}.

\paragraph{Minimax optimization.}
Finite-sum saddle-point methods, matrix games, and variance-reduced
extragradient \cite{saddle2016,matrixgames2019,svre2019,catalystminimax2020}; recursive, hybrid, adaptive, and decentralized
nonconvex minimax methods \cite{sreda2020,hybridminimax2020,adgda2023,dream2024}; a variational-inequality treatment \cite{vi2022}.


\end{document}